\documentclass[a4paper,11pt]{amsart}
\usepackage[left=3.6cm, top=3cm, bottom=3.5cm, right=3.6cm]{geometry}
\usepackage{amsthm,amsmath}
\usepackage{amssymb}
\usepackage{mathtools}
\usepackage{enumitem}
\usepackage{graphicx}
\usepackage{microtype}
\usepackage[T1]{fontenc}
\usepackage{booktabs}
\usepackage{float}
\usepackage{lmodern}

\newtheorem{thm}{Theorem}
\newtheorem{pro}[thm]{Proposition}

\newtheorem{lem}[thm]{Lemma}

\theoremstyle{definition}
\newtheorem{defn}[thm]{Definition}
\newtheorem{rem}[thm]{Remark}
\newtheorem{exm}[thm]{Example}

\numberwithin{thm}{section}
\numberwithin{equation}{section}
\numberwithin{figure}{section}
\numberwithin{table}{section}

\makeatletter
\def\thmhead@plain#1#2#3{%
  \thmname{#1}\thmnumber{\@ifnotempty{#1}{ }\@upn{#2}}%
  \thmnote{ {\the\thm@notefont#3}}}
\let\thmhead\thmhead@plain
\makeatother

\def\del{\partial}

\def\C{\mathbb{C}}
\def\HH{\mathbb{H}}
\def\N{\mathbb{N}}

\def\Q{\mathbb{Q}}
\def\R{\mathbb{R}}
\def\Z{\mathbb{Z}}
\def\CP{\mathbb{C}\text{P}}
\def\SL{\text{SL}}

\def\sC{\mathcal{C}}
\def\sD{\mathcal{D}}
\def\sE{\mathcal{E}}

\def\sL{\mathcal{L}}

\def\sR{\mathcal{R}}

\def\sW{\mathcal{W}}

\def\sZ{\mathcal{Z}}

\def\fX{\mathfrak{X}}
\def\fg{\mathfrak{g}}
\newcommand{\betti}{{\operatorname{b}}}
\def\fg{\mathfrak{g}}
\newcommand{\kod}{{\operatorname{kod}}}
\def\JE{\mbox{$J$-Engel }}

\def\bra{\langle}
\def\ket{\rangle}

\newcommand{\Matrix}[1]
{
\left(
\begin{matrix}
#1
\end{matrix}
\right)
}

\newcommand{\su}[1]
{
_{_{#1}}
}

\usepackage{nicefrac}
\makeatletter
\newsavebox\myboxA
\newsavebox\myboxB
\newlength\mylenA
\newcommand*\xoverline[2][0.75]{%
    \sbox{\myboxA}{$\m@th#2$}%
    \setbox\myboxB\null
    \ht\myboxB=\ht\myboxA%
    \dp\myboxB=\dp\myboxA%
    \wd\myboxB=#1\wd\myboxA
    \sbox\myboxB{$\m@th\overline{\copy\myboxB}$}
    \setlength\mylenA{\the\wd\myboxA}
    \addtolength\mylenA{-\the\wd\myboxB}%
    \ifdim\wd\myboxB<\wd\myboxA%
       \rlap{\hskip 0.5\mylenA\usebox\myboxB}{\usebox\myboxA}%
    \else
        \hskip -0.5\mylenA\rlap{\usebox\myboxA}{\hskip 0.5\mylenA\usebox\myboxB}%
    \fi}
\makeatother

\title{Engel structures on complex surfaces}

\author{N. Pia}
\address{Mathematisches Institut, LMU M\"unchen, Theresienstr. 39, 80333 M\"unchen, Germany}
\email{nicola.pia@hotmail.it}
\author{G. Placini}
\address{Mathematisches Institut, LMU M\"unchen, Theresienstr. 39, 80333 M\"unchen, Germany}
\email{placini@math.lmu.de}

\subjclass[2010]{Primary 53C56; Secondary 32Q99, 53D99} 
\keywords{Engel structures; Complex surfaces.}

\begin{document}
%
%
\begin{abstract}
We classify complex surfaces $(M,\,J)$ admitting Engel structures $\mathcal{D}$ which are complex line bundles.
Namely we prove that this happens if and only if $(M,\,J)$ has trivial Chern classes.
We construct examples of such Engel structures by adapting a construction due to Geiges~\cite{geigesEngel}.
We also study associated Engel defining forms and define a unique splitting of $TM$ associated with $\mathcal{D}$ \mbox{$J$-Engel}.
\end{abstract}

\maketitle

%
%
\section{Introduction}
An Engel structure $\sD$ is a maximally non-integrable $2$-plane field on a $4$-manifold $M$, i.e. $[\sD,\sD]=\sE$ and $[\sD,\sE]=TM$.
These objects were discovered a long time ago \cite{cartan, engel}, and recent developments have sparked new interest in the field \cite{cp3, overtwistedEngel, vogelExistence}.

Line fields, contact structures, even contact structures, and Engel structures are the only topologically stable families of distributions in the sense of Cartan \cite{cartan}.
This means in particular that they admit Darboux-type theorems, which implies that they do not have local invariants.
Engel structures are an exceptional family in this list since they only exist on $4$-dimensional (virtually) parallelizable manifolds.

This paper concerns the interplay between Engel structures and complex structures.
For a given almost complex $4$-manifold $(M,\,J)$, a \JE structure is an Engel structure $\sD$ such that $J\sD=\sD$.
These structures have already been studied in the case where $J$ is integrable in \cite{zhao}, where Zhao studies the case of homogeoneous \JE structures and classifies the structure constants in this case.

In Zhao's work \JE structures appear under the name of \emph{complex Engel structures}.
We prefer to use the former name in order not to create confusion with holomorphic Engel structures \cite{coelhopia, presasconde}, which are the analogue of Engel structures in the holomorphic category.

Notice that if $\sD$ is \JE then we can find a vector field $W$ tangent $\sD$ such that $\big\{ W,\,JW,\,[W,JW],\,J[W,JW]\big\}$ is a framing for $TM$, which implies that the Chern classes of $(M,\,J)$ vanish.
The main result in this paper, answering a question in \cite{zhao}, is the following.
\begin{thm}\label{THM_Main}
  A complex surface $(M,J)$ admits a \mbox{$J$-Engel} structure if and only if $c_1(M)=0$ and $c_2(M)=0$.
\end{thm}
The strategy of the proof is to use Enriques-Kodaira classification to rule out complex surfaces with inadmissible Chern classes.
It turns out that all complex surfaces with trivial Chern classes are (orientable) mapping tori.
This allows us to adapt a construction of Engel structures on mapping tori due to Geiges~\cite{geigesEngel}, hence proving Theorem~\ref{THM_Main}.

\begin{rem}
A statement analogous to Theorem~\ref{THM_Main} holds true when replacing \JE structure by totally real orientable Engel structure, cf. Remark~\ref{totreal}.
\end{rem}

Towards the end of the paper we analyse some properties of defining forms associated with a \JE structure and present a list of interesting examples.
Two $1$-forms $\alpha$ and $\beta$ are said to be Engel defining forms for a given Engel structure $\sD$ if $\sD=\ker\alpha\cap\ker\beta$ and $\sE=[\sD,\sD]=\ker\alpha$.
A pair of defining forms determines a complementary distribution $\sR=\bra T,\,R\ket$ called the Reeb distribution (see \cite{nicola}).
The conformal class of $\alpha$ is uniquely determined by $\sD$, whereas in general the conformal class of $\beta$ is not.
If $\sD$ is a \JE structure we do have a natural choice for this conformal class, namely $\beta=\alpha\circ J$.
This in turn defines the line bundle $\sZ=\bra R\ket$ uniquely, and hence we get a splitting of the tangent bundle of $M$ into a sum of subline bundles
\[
 TM=\sW\oplus J\sW\oplus J\sZ\oplus \sZ
\]
called the \JE framing.
We study the case where the flow of the vector field $R$ acts by \JE isomorphisms and list some geometric examples.

\subsection{Structure of the Paper}
In Section~\ref{SEC_BasicFacts} we recall some basic facts from the theory of Engel structures.
Section~\ref{SEC_Proof} is dedicated to the proof of Theorem~\ref{THM_Main}.
There we give the details of the classification of surfaces with trivial Chern classes and we construct \JE structures on complex mapping tori.
In Section~\ref{SEC_Forms} we study properties of defining forms associated with a \JE structure, define the \JE framing, and study the case where the flow of $R$ acts via \JE isomorphisms.
Finally Section~\ref{SEC_Examples} contains interesting constructions of \JE structures with a focus on Thurston geometries.

\subsection{Acknowledgements}
We would like thank our advisor Prof. Kotschick for the useful discussions and for pointing out to us the results in \cite{Li, teleman, wall}.
Moreover we thank Rui Coelho for the helpful and clarifying conversations.

%
%
\section{Engel structures}\label{SEC_BasicFacts}
In what follows all manifolds are assumed to be closed and smooth, and all distributions are assumed to be smooth, if not otherwise stated.

An \emph{even contact structure} $\sE$ is a maximally non-integrable hyperplane distribution on an even dimensional manifold.
Otherwise said, if $\dim (M)=2n+2$ then locally $\sE$ is the kernel of a 1-form $\sE=\ker\alpha$ satisfying $\alpha\wedge d\alpha^{2n}\ne0$.
For dimensional reasons if $\sE$ is even contact then there exists a unique line field $\sW$ such that $\sW\subset\sE$ and $[\sW,\sE]\subset\sE$.
We call $\sW$ the \emph{characteristic foliation} of $\sE$.

An \emph{Engel structure} $\sD$ is a smooth $2$-plane field on a smooth $4$-manifold such that $\sE:=[\sD,\sD]$ is an even contact structure.
One can see that the characteristic foliation $\sW$ of $\sE=[\sD,\sD]$ satisfies $\sW\subset\sD$.
The flag of distributions $\sW\subset\sD\subset\sE$ is called the \emph{Engel flag of} $\sD$.
The existence of this flag gives strong constraints on the topology of the manifold $M$.
If $(M,\sD)$ is an Engel structure and $\sW\subset\sD\subset\sE$ is its associated flag then we have canonical isomorphisms
\begin{equation}\label{EQN_IsoOfDetBundles}
	\det(\sE/\sW)\cong\det(TM/\sE)\qquad\textup{and}\qquad\det(\sE/\sD)\cong\det(\sD).
\end{equation}
In particular $\sE$ is orientable and $M$ is orientable if and only if $\sW$ is trivial.
This implies that if $M$ admits an Engel structure, then it admits a parallelizable $4$-cover.
Notice that if $M$ is orientable and $\sD$ is orientable then we can construct a framing $TM=\bra W,\,X,\,Y,\,Z\ket$ such that
$\sW=\bra W\ket$, $\sD=\bra W,\,X\ket$ and $\sE=\bra W,\,X,\,Y\ket$.

We recall the following characterisation of parallelizable $4$-manifolds
\begin{thm}[\cite{hiho}]
 An orientable $4$-manifold is parallelizable if and only if its Euler characteristic $\chi(M)$, second Stiefel-Whitney class $w_2(M)$, and signature $\sigma(M)$ vanish.
\end{thm}

It was an open question for a long time whether all parallelizable manifold admit Engel structures.
This problem was solved positively for the first time in \cite{vogelExistence}.
The later works \cite{cp3, overtwistedEngel} established an existence h-principle and constructed a flexible (in the sense of Gromov) family of Engel structures.

\begin{defn}
 Let $(M^4,J)$ be an almost complex manifold, a $J$-Engel structure $\sD$ is an Engel structure such that $J\sD=\sD$.
\end{defn}
In this paper we are interested in the study of Engel structures which are complex line fields on complex manifolds.
Throughout the paper $(M,J)$ will denote a closed four-manifold $M$ equipped with a complex structure $J$ if not otherwise specified.

The first example of \mbox{$J$-Engel} structure on a complex surface was constructed by Bryant during the workshop \emph{Engel structures} in San Jos\'e in 2017.
      \begin{exm}[(R.~Bryant)]
      Let $z,w$ be holomorphic coordinates on $\mathbb{C}^2$ and define 
      \[
      \omega=e^{i(z+\bar z)}dw-idz
      \]
      Let $\mathbb{Z}^4$ act on $\mathbb{C}^4$ via
      \[
      (k_1,k_2,k_3,k_4)(z,w)=(z+\pi k_1+ik_2,w+k_3+ik_4)
      \] 
      Then $\omega$ passes to the quotient $M=\C^2/\Z^4=T^4$, defining on it an Engel structure $\sD=\ker\omega$.
      \end{exm}

%
%
\section{Proof of Theorem~\ref{THM_Main}}\label{SEC_Proof}

We now investigate the topological constraints on $(M,J)$ given by the existence of a \mbox{$J$-Engel} structure.

Suppose that an almost complex manifold $(M,J)$ admits a \mbox{$J$-Engel} structure.
Since $M$ is almost complex, it is orientable, hence Equation~\eqref{EQN_IsoOfDetBundles} implies that $\sW$ is trivial as a bundle.
Fix a non-vanishing section $W\in\Gamma\sW$.
Since $\sD$ is \mbox{$J$-Engel}, we have that $X=JW$ is tangent to $\sD$ and $\sD=\bra W,\,X\ket$.
Maximal non-integrability implies that $Y=[W,X]$ is nowhere tangent to $\sD$, so that we get the complex framing $TM=\bra W,\,JW,\,Y,\,JY\ket$. This shows that the Chern classes of $(M,J)$ are trivial and in particular concludes the proof of the necessity in Theorem~\ref{THM_Main}.

We now prove the converse. In order to list all complex surfaces with trivial Chern classes we make use of Enriques--Kodaira classification which we recall here:

\begin{table}[H]
\begin{tabular}{@{}llccl@{}}
\toprule
$\kod(M)$ \hspace{1mm}	& Class of $M$  \hspace{1cm}	& $c_1^2(M)$ \hspace{5mm}	 & $c_2(M)$		 	\\ \midrule
 		    & minimal rational surfaces  		       & $8,9$	   & $4,3$       \\
 negative	& class VII minimal surfaces 	           & $\leq0$   & $\geq0$     \\
			& minimal ruled surfaces of genus $g\geq1$ & $8(1-g)$  & $4(1-g)$    \\ \midrule
 			& Enriques surfaces                        & $0$       & $12$        \\
 	        & hyperelliptic surfaces                   & $0$       & $0$         \\
	$0$		& Kodaira surfaces		                   & $0$       & $0$         \\
 			& $K3$-surfaces                            & $0$	   & $24$        \\
        	& complex tori	                           & $0$	   & $0$         \\ \midrule
	$1$		& minimal properly elliptic surfaces    	           & $0$       & $\geq0$     \\ \midrule
    $2$		& surfaces of general type                 & $>0$      & $>0$        \\ \bottomrule
\end{tabular}
\caption{\label{classification}Enriques--Kodaira classification.}
\end{table}
In the above table $\kod(M)$ denotes the Kodaira dimension of the surface $M$.
We can now prove the following
\begin{lem}\label{LEM_TopologyOfSurfacesWithTrivialChern}
  If a complex surface $(M,J)$ has $c_1(M)=0$ and $c_2(M)=0$ then $M$ is minimal and belongs to one of the following families
  \begin{itemize}
   \item Inoue surfaces;
   \item Hopf surfaces;
   \item Hyperelliptic surfaces;
   \item Kodaira surfaces;
   \item Complex tori;
   \item Non-K\"ahler surfaces of Kodaira dimension 1.
  \end{itemize}

  \begin{proof}
  We refer to Table~\ref{classification} for the invariants of the different classes of surfaces.
  
  We can immediately exclude rational surfaces, Enriques surfaces, $K3$-surfaces, and surfaces of general type, since the Euler characteristic increases under blow-ups and their minimal models have positive Euler characteristic. 

  Ruled surfaces are birationally equivalent to $\CP^1\times\Sigma_g$, where $\Sigma_g$ denotes a curve of genus $g$.
  The hypothesis imply that the signature of $M$ must vanish, hence if $M$ is ruled it must be minimal.
  In this case though $c_1(M)\ne0$.

  In all other classes the only surfaces with $c_2(M)=0$ are minimal.
  Notice that class ${\rm VII}$ surfaces $M$ have $\betti_1(M)=1$, hence $c_2(M)=0$ implies $\betti_2(M)=0$ in this case.
  It was proven that only Hopf surfaces and Inoue surfaces satisfy this condition, see~\cite{Li,teleman}.
  Finally a K\"ahler elliptic surface with trivial first Chern class must have Kodaira dimension $0$.
  Indeed if the canonical bundle has a non-vanishing section, then this defines a hypersurface dual to the first Chern class, i.e., nullhomologous, contradicting the fact that the K\"ahler class evaluates non-trivially on it.
  \end{proof}
\end{lem}

All surfaces in the classes listed in the above lemma have the structure of a mapping torus, i.e. they are fibre bundles over $S^1$.
This is clear for hyperelliptic surfaces and tori.
The same is true for Inoue surfaces, cf.~\cite{inoue}, while the statement for Hopf surfaces was proven by Kato in \cite{kato, katoErratum}.
The remaining classes, i.e. Kodaira surfaces and non-K\"ahler surfaces of Kodaira dimension $1$, admit a Vaisman metric (see \cite{belgun} for the classification).
Ornea and Verbitsky~\cite{ornea} proved that all such manifolds are mapping tori with Sasakian fibers.

In order to prove Theorem~\ref{THM_Main} it suffices to show that a complex surface with trivial first Chern class which fibers over $S^1$ admits a \JE structure.
This is the content of Lemma~\ref{LEM_JEngelGeiges}, which we state in the more general setting of almost complex manifolds.
\begin{rem}
 Notice that all surfaces in the families listed in Lemma~\ref{LEM_TopologyOfSurfacesWithTrivialChern} have trivial first Chern class, except possibly for Hopf surfaces.
 In this case indeed it is unclear if all secondary Hopf surfaces satisfy $c_1(M)=0$, notice that this is true for primary ones since they are diffeomorphic to $S^1\times S^3$ (see \cite{kato}). 
\end{rem}

\begin{lem}\label{LEM_JEngelGeiges}
 Let $(M,J)$ be an almost complex $4$-manifold such that $c_1(M)=0$.
 Moreover suppose that $M$ is diffeomorphic to an orientable $3$-manifold bundle over $S^1$, then $M$ admits a \mbox{$J$-Engel} structure.
 \begin{proof}
    Let $\pi\colon M\longrightarrow  S^1$ be the bundle projection with (oriented) fibre $N$ and denote by $f\colon N\longrightarrow  N$ the monodromy of the bundle.
    This means that $M$ is the suspension of the (orientation preserving) diffeomorphism $f\colon N\longrightarrow  N$.
    The vector field $\partial_t$ on $S^1$ induces a nowhere-vanishing vector field $V$ on $M$ which satisfies $\pi_*V=\partial_t$.
    This implies that $c_2(M)=0$.
    Since $c_1(M)=0$ we have a framing of the form $TM=\bra V,\,JV,\,X,\,JX\ket$ for some vector field $X\in\fX(M)$.
    Let $n\in\N$ and consider the plane field $\sD_n=\bra A_n,\,JA_n\ket$ where
    \[
     A_n :=V+\frac{1}{n}\sin(n^2\pi^*t)\ X-\frac{1}{n}\cos(n^2\pi^*t)\ JX.
    \]
    Using $\sL_V\pi^*t=1$ and $a:=\sL_{JV}\pi^*t$ we have
    \begin{align*}
     JA_n&=JV+\frac{1}{n}\cos(n^2\pi^*t)\ X+\frac{1}{n}\sin(n^2\pi^*t)\ JX\\
     \frac{1}{n}[A_n,JA_n]&=-\Big(\sin(n^2\pi^*t)+a\cos(n^2\pi^*t)\Big)X\\
     &\quad +\Big(\cos(n^2\pi^*t)-a\sin(n^2\pi^*t)\Big)JX+O(n^{-1})\\
     \frac{1}{n^2}\left[A_n,\frac{1}{n}[A_n,JA_n]\right]&=-\Big(\cos(n^2\pi^*t)-a\sin(n^2\pi^*t)\Big)X\\
     &\quad -\Big(\sin(n^2\pi^*t)+a\cos(n^2\pi^*t)\Big)JX+O(n^{-1}).
    \end{align*}
    Hence if $n$ is big enough $\sD_n$ is a \mbox{$J$-Engel} structure. 
 \end{proof}
\end{lem}
This completes the proof of Theorem~\ref{THM_Main}.
In \cite{geigesEngel} Geiges exhibited an Engel structure on parallelizable manifolds $M$ which are suspensions of a diffeomorphism $f\colon N\longrightarrow  N$ of a $3$-manifold $N$.
Lemma~\ref{LEM_JEngelGeiges} adapts this construction to the \mbox{$J$-Engel} case.

\begin{rem} \label{totreal}
Instead of $\sD$ being complex one can ask the Engel distribution $\sD$ to be totally real, i.e. $J\sD\cap\sD=0$.
If $\sD$ is totally real and orientable the analogue of Theorem~\ref{THM_Main} still holds.
In fact such an Engel structure $\sD=\bra W,\, X\ket$ is trivial as a bundle, so that $TM=\bra W,\,JW,\,X,\,JX\ket$ proving that the same restrictions on the Chern classes apply.
On the other hand, with the notation of Lemma~\ref{LEM_JEngelGeiges}, the distribution $$\sD=\bra V,\,JV+\frac{1}{n}\cos(n^2\pi^*t)\ X+\frac{1}{n}\sin(n^2\pi^*t)\ JX\ket$$ gives a totally real Engel structure on all complex surfaces $M$ with $c_1(M)=0$ and $c_2(M)=0$.
\end{rem}

%
\section{Engel defining forms}\label{SEC_Forms}
The first author has studied the properties of particular $1$-forms $\alpha$ and $\beta$ such that the intersection of their kernels is an Engel structure (see \cite{nicola}), we now consider these objects in the context of \JE structures.
The results that follow are true also in the case where $J$ is non-integrable.

Let $\sD$ be an Engel structure. If two $1$-forms $\alpha$ and $\beta$ satisfy $\sD=\ker\alpha\cap\ker\beta$ and $\sE=\ker\alpha$, we say that $\alpha$ and $\beta$ are \emph{Engel defining forms for} $\sD$.
This happens if and only if
\[
	\alpha\wedge d\alpha\ne0,\quad \alpha\wedge\beta\wedge d\beta\ne0,\quad\textup{and}\quad\alpha\wedge d\alpha\wedge\beta=0.
\]
A pair of defining forms determines a distribution $\sR=\bra T,\, R\ket$ transverse to $\sD$ via
\[
\begin{matrix}
  i_T(\alpha\wedge d\beta)=0,&\beta(T)=1,&\alpha(T)=0, \\
  i_R(\beta\wedge d\beta)=0,&\beta(R)=0,&\alpha(R)=1.
\end{matrix}
\]
This is called the \emph{Reeb distribution associated with $\alpha$ and $\beta$}.
An Engel structure $\sD$ on an orientable manifold $M$ admits Engel defining forms if and only if it is orientable.
This is the case for \JE structures.

A direct calculation proves that if $\alpha$ is a defining form for $\sE$, then we can get a pair of Engel defining forms for $\sD$ by setting $\beta=\alpha\circ J$.
Any other choice of $\alpha$ is of the form $\lambda\,\alpha$ for some nowhere vanishing function $\lambda$, hence the same is true for $\beta=\lambda\,\alpha\circ J$.
This simple observation implies that a \JE structure gives a preferred splitting of the tangent bundle, as the following result ensures.
\begin{pro}\label{PROP_JESplitting}
 Suppose that $\sD$ is a \JE structure on an almost complex $4$-manifold $(M,J)$.
 Let $\sW$ be the characteristic foliation of $\sD$, let $\alpha$ be a defining form for $\sE$, $\beta=\alpha\circ J$, and denote by $\sZ=\bra R\ket$.
 The splitting of the tangent bundle
 \begin{equation}\label{EQN_JESplitting}
  TM=\sW\oplus J\sW\oplus J\sZ\oplus\sZ,
 \end{equation}
 does not depend on the choice of $\alpha$.
 We call it the \JE splitting.
 
 \begin{proof}
  All other possible choices of $\alpha$ are of the form $\tilde\alpha=\lambda\,\alpha$ for $\lambda$ nowhere-vanishing function.
  This implies that $\tilde\beta=\lambda\,\beta$ and $\tilde\beta\wedge d\tilde\beta=\lambda^2\,\beta\wedge d\beta$, so that $\lambda\tilde R=R$ concluding the proof.
 \end{proof}

\end{pro}
The previous result implies that any isomorphism $\phi:M\to M$ which preserves both the Engel structure $\sD$ and the complex structure $J$, must also preserve the associated \JE splitting.
We say that a vector field $Z$ is \JE if its flow preserves both $\sD$ and $J$.

An interesting question in the field is whether every Engel structure admits a $1$-parameter family of symmetries, i.e. a vector field $Z$ whose flow preserves the Engel structure, also called Engel vector field \cite{montgomeryPaper, nicola}.
The following result gives a necessary condition in the \JE setting.
\begin{lem}\label{LEM_TransverseEngelVF}
 Suppose that $\sD$ is a \JE structure on an almost complex $4$-manifold $(M,J)$.
 Let $Z$ be an Engel vector field transverse to $\sE$ and such that $JZ\in\Gamma\sE$, then $Z\in\Gamma\sZ$.
 
 In particular there exists $\alpha$ and $\beta=\alpha\circ J$ Engel defining forms such that $Z=R$.
 \begin{proof}
  Fix any defining form $\alpha$ for $\sE$ and consider $\beta=\alpha\circ J$, so that the Reeb distribution $\sR=\bra T,\, R\ket$ is well defined.
  To prove the claim it suffices to show that $i_Z(\beta\wedge d\beta)=0$, that is $d\beta(Z,A)=0$ for $A\in\ker\beta$.
  In fact, since the formula is verified by definition for $A=R$, it is enough to prove it for $A\in\sD$.
  Since $\beta(Z)=\alpha(JZ)=0$ by hypothesis, Cartan formula gives
  \[
  d\beta(Z,A)=-\beta([Z,A])=0,
  \]
  where we used the hypothesis $\sL_Z\sD\subset\sD$.
  
  The proof of the second claim follows by taking the Engel defining forms
  \[
   \tilde\alpha=\frac{1}{\alpha(Z)}\,\alpha\quad\textup{and}\quad\tilde\beta=\tilde\alpha\circ J.
  \]

 \end{proof}
\end{lem}

We now turn back to the case where $J$ is integrable, and describe the interplay between $J$ and the Reeb distribution.
We fix a framing $\sD=\bra W,\,X=JW\ket$ and Engel defining forms $\alpha$ and $\beta$. In the following we use the notation from \cite[Section~3]{nicola}, that is $c_{_{WX}}=\beta([W,X])$, $d_{_{XT}}=\alpha([X,T])$, $d_{_{WR}}=\alpha([W,R])$, and $d_{_{XR}}=\alpha([X,R])$.
\begin{lem}\label{LEM_JofReeb}
 Suppose that $\sD$ is a \JE structure on a complex surface $(M,J)$, let $\alpha$ be a defining form for $\sE$, and $\beta=\alpha\circ J$.
 Then the Reeb distribution $\sR=\bra T,\, R\ket$ satisfies
 \begin{align*}
  JT&=R+\frac{d_{_{WR}}+d_{_{XT}}}{c_{_{WX}}}\,W+\frac{d_{_{XR}}}{c_{_{WX}}}\,JW\\
  JR&=-T+\frac{d_{_{XR}}}{c_{_{WX}}}\,W-\frac{d_{_{WR}}+d_{_{XT}}}{c_{_{WX}}}\,JW.
 \end{align*}

 \begin{proof}
  Suppose that $JT=a\,W+b\,X+c\,T+d\,R$ for $a,b,c,d\in\sC^\infty(M)$, we have
  \[
   c=\beta(JT)=-\alpha(T)=0\quad\textup{and}\quad d=\alpha(JT)=\beta(T)=1.
  \]
  Now since $J$ is integrable we have
  \[
   [W,JT]=J[W,T]-J[JW,JT]-[JW,T],
  \]
  which in turn yields
  \begin{align*}
   \alpha([W,R])&=\alpha([W,JT-a\,W-b\,X])=\alpha([W,JT])\\
   &=\alpha(J[W,T])-\alpha(J[JW,JT])-\alpha([JW,T])\\
   &=\beta([W,T])-\beta([JW,R+a\,W+b\,X])-d_{_{XT}}\\
   &=a\,\beta([W,JW])-1=a\,c_{_{WX}}-d_{_{XT}}.
  \end{align*}
  The formula for $b$ is obtained via a similar calculation for $\alpha([X,R])$, and the formula for $R$ is a consequence of $J^2=-id$.
 \end{proof}
\end{lem}

\begin{rem}
Notice that $d\alpha^2=-2d_{_{WR}}\,\alpha\wedge\beta\wedge d\beta$, this means that, since $M$ is closed, there are points where the function $d_{_{WR}}$ vanishes.
This implies that the Reeb distribution associated with the forms in the previous lemma is never $J$-invariant.
 \end{rem}

We now study the very special case where there exists a \JE vector field $Z$ transverse to $\sE$ and such that $JZ\in\Gamma\sE$.
In view of Lemma~\ref{LEM_TransverseEngelVF}, this is the same as saying that we have a defining form $\alpha$ such that the flow of $R$ acts by \JE isomorphisms. This particular instance provides a connection between \JE structures and K-Engel structures.

A K-Engel structure on $M$ is a triple $(\sD,\,g,\,Z)$ where $\sD$ is an Engel structure, $g$ a Riemannian metric, and $Z$ a Killing Engel vector field which is orthogonal to $\sE$.
In order to prove that $\sD$ admits a K-Engel structure it suffices to exhibit Engel defining forms $\alpha$ and $\beta$ and a framing $\sD=\bra W,\,X\ket$ such that $R$ commutes with $W$, $X$, and $T$ (see \cite[Proposition 7.5]{nicola}).
\begin{pro}\label{PROP_JEKE}
 Suppose that $\sD$ is a \JE structure on a complex surface $(M,J)$, and let $Z$ be a \JE vector field transverse to $\sE$ and such that $JZ\in\Gamma\sE$, then $\sD$ admits a K-Engel structure.
 
 \begin{proof}
  Using Lemma~\ref{LEM_TransverseEngelVF} we can suppose that $\alpha$ and $\beta=\alpha\circ J$ are Engel defining forms such that $Z=R$.
  It suffices to find a section $W$ of the characteristic foliation such that $R$ commutes with $W$, $X=JW$ and $T$.
  Our hypotheses ensure that $d\su{WR}=0=d\su{XR}$ so that
  \begin{align*}
   [W,R]&=a\su{WR}W\\
   [X,R]&=a\su{XR}W+b\su{XR}X\\
   [T,R]&=a\su{TR}W+b\su{TR}X+c\su{TR}T.
  \end{align*}
  Since $X=JW$ and $\sL_RJ=0$ we conclude that $a\su{XR}=0$ and $b\su{XR}=a\su{WR}$ .
  Now Lemma \ref{LEM_JofReeb} ensures that
  \[
   JR=-T-\frac{d_{_{XT}}}{c_{_{WX}}}\,JW=:-T-aX.
  \]
  A direct calculation shows
  \[
   0=[R,JR]=[T,R]+[aX,R]=a\su{TR}W+b\su{TR}X+c\su{TR}T+(\sL_Ra)\,X+a\,a\su{WR}\,X,
  \]
  in particular $c\su{TR}=0$ implying $d\beta^2=0$, so that $\sR$ is a foliation and hence $a\su{TR}=0=b\su{TR}$ (for more details see \cite[Section 3]{nicola}).
  Moreover we have
  \[
   a\su{WR}=-\sL_R(\ln a)
  \]
  which implies that we can rescale $W$ so that $a\su{WR}=0$, hence proving the statement.
 \end{proof}
\end{pro}

%
\section{Examples of \mbox{$J$-Engel} structures}\label{SEC_Examples}
This section is dedicated to constructing explicit examples of \mbox{$J$-Engel} structures on complex surfaces.

Some classes admit geometric structures that can be used to produce such examples.
More precisely $M$ admits a geometric structure if it is modelled on a simply connected manifold $X$ with a transitive action of a Lie group $G$ and a $G$-invariant metric.
A geometric complex structure on $M$ is a $G$-invariant complex structure.
In our case, since $M$ is compact, $X=\widetilde M$ there is a lattice $\Gamma$ of $G$ such that $M\cong X/\Gamma$.

The following theorem of Wall classifies geometric structures on (not necessarily properly) elliptic surfaces.

\begin{thm}[\cite{wall}]\label{THM_geometric}
An elliptic surface $M$ without singular fibres has a geometric structure if and only if its base is a good orbifold\footnote{
This means that it admits a finite orbifold cover with no cone points (see Section 7 in \cite{wall} for more details).}.
The geometric structure is determined as follows

\begin{table}[H]
\begin{tabular}{@{}lcccl@{}}
\toprule
$\kod(M)$ \hspace{2cm}	& negative   \hspace{1cm}	& $0$ 	\hspace{1cm} & $1$		 	\\ \midrule
 $\betti_1(M)$ even	& $\C\times\CP^1$  \hspace{1cm}  & $\C^2$	 \hspace{1cm}  & $\C\times\HH$    \\
 $\betti_1(M)$ odd	& $S^3\times \R$   \hspace{1cm}  & $Nil_3\times \R$\hspace{1cm} & $\widetilde{\SL(2,\R)}\times\R$        \\ \bottomrule
\end{tabular}
\end{table}
\end{thm}

Moreover the following result classifies the ones which admit solv geometries.
\begin{thm}[\cite{hasegawa}]\label{THM_hasegawa}
 A complex surface is diffeomorphic to a $4$-dimensional solvmanifold if and only if it is one of the following surfaces: complex torus, hyperelliptic surface, Inoue surface of type $S_0$, primary Kodaira surface, secondary Kodaira surface, Inoue surface of type $S_\pm$.
 And every complex structure on each of these complex surfaces (considered as solvmanifolds) is left-invariant.
\end{thm}
Moreover, in \cite[Section~5]{hasegawa}, Hasegawa gives an explicit construction of the complex structure on these solvmanifolds.
Namely consider a $4$-dimensional simply connected solvable Lie group $G$, so that $G/\Gamma$ is a solvmanifold for $\Gamma$ cocompact lattice.
Let $\fg$ be the Lie algebra of $G$, fix a basis $\{X_1,\,X_2,\,X_3,\,X_4\}$, and construct an almost complex structure $J$ by defining
\[
  JX_1 = X_2,\quad JX_3 = X_4 \ .
\]
The following list gives the left-invariant complex structures in Theorem~\ref{THM_hasegawa} 
\begin{enumerate}
  \item \emph{Complex Tori}: $G=\R^4$ and all brackets vanish.
  \item \emph{Hyperelliptic surfaces}: all brackets vanish except for
  \[
   [X_1,X_4] = X_2,\quad [X_2,X_4] = - X_1.
  \]
  \item \emph{Primary Kodaira surfaces}: $G=Nil_3\times\R$ and all brackets vanish except for
  \[
   [X_1,X_2] = - X_3.
  \]
  \item \emph{Secondary Kodaira surfaces}: $G$ is the maximal connected isometry group of $Nil_3$ and all brackets vanish except for
  \[
   [X_1,X_2] = - X_3,\quad [X_1,X_4] = X_2,\quad [X_2,X_4] = - X_1.
  \]
  \item \emph{Inoue surfaces of type $S_0$}: $G=Sol_0^4$ and all brackets vanish except for
  \[
   [X_1,X_4] = -a X_1 + b X_2,\quad [X_2,X_4] = - b X_1 - a X_2,\quad [X_3,X_4] = 2 a X_3,
  \]
  where $a,b\in\R^*$.
  \item \emph{Inoue surfaces of type $S_+$ and $S_{-}$}: $G=Sol_1^4$ and all brackets vanish except for
  \[
   [X_2,X_3] = - X_1,\quad [X_2,X_4] = - X_2,\quad [X_3,X_4] = X_3.
  \]
  In this case there is a family of almost complex structures given by
  \[
    JX_1 = X_2,\quad JX_3 = X_4 - q X_2,
  \]
  for $q\in\R$.
\end{enumerate}
One can verify that the above formulae define integrable almost complex structures.

In the remainder of this section we give examples of \JE structures in each family appearing in Lemma~\ref{LEM_TopologyOfSurfacesWithTrivialChern} and point out which of these examples are K-Engel.

\subsection{Inoue surfaces}\label{SEC_Inoue}
Let us consider first Inoue surfaces of type $S_0$ as solvmanifolds with the Lie algebra structure given above.
A left-invariant \mbox{$J$-Engel} structure is given by $\sD=\bra A ,\, JA\ket$ where $A=X_1+X_4$. Indeed a simple computation yields
\begin{align*}
[A,\, JA]&=bX_1+aX_2+2aX_3\, , \\
[A,\, [A,\, JA]]&=2abX_1+(a^2-b^2)X_2-4a^2X_3\, ,
\end{align*}
and these four vectors span the Lie algebra.
Notice that \cite[Section 11]{nicola} implies that in this case we do not have any geometric K-Engel defining forms.

For Inoue surfaces of type $S_{\pm}$ one considers the left-invariant complex plane field $\sD=\bra A ,\, JA\ket$ where $A=X_1+X_4$.
The Lie algebra structure then gives
\begin{align*}
[A,\, JA]&=X_2+X_3\,\, \text{and} \\
[JA,\, [A,\, JA]]&=-2X_1\, ,
\end{align*}
which shows that $\sD$ is a \JE structure.

\subsection{Hopf surfaces}\label{SEC_Hopf}
Some Hopf surfaces are modelled on $G=S^3\times\R$.
We can fix a basis $\{X_1,\,X_2,\,X_3,\,X_4\}$ for the Lie algebra $\fg$ of $G$ for which the only non-zero Lie brackets are
\[
 [X_1,X_2]=X_3\qquad[X_2,X_3]=X_1\qquad[X_3,X_1]=X_2.
\]
and the complex structure $J$ is given by $JX_1=X_2$ and $JX_3=X_4$.
We denote by $\{a_1,\,a_2,\,a_3,\,a_4\}$ the dual basis on the dual Lie algebra.

We can define a left-invariant complex plane field on $G$ by $\sD=\bra A,\,JA\ket$ with $A=X_1+X_3$.
By computing
\begin{align*}
[A,\, JA]&=X_3-X_1\,\, \text{and} \\
[A,\, [A,\, JA]]&=X_2\, 
\end{align*}
one sees that $\sD$ defines a \JE structure.
In this case $R=X_4$ is a \JE vector field such that $\alpha(JR)=0$, so $\sD$ admits a K-Engel structure, in fact the Engel defining forms are given by $\alpha=a_4-a_2$ and $\beta=\alpha\circ J=a_3-a_1$.

\subsection{Hyperelliptic surfaces}
Any hyperelliptic surface $M$ is the quotient of the product of two elliptic curves $T^2_\Lambda\times T^2_{\tilde\Lambda}$ by the action of a finite group $G$ (see for instance \cite{barth}).
More explicitly we take coordinates $z_j=x_j+iy_j$ for $j=1,2$ on $\C^2$, and we denote by $\omega$ the primitive third root of the identity.
The admissible finite groups $G$ and their actions were classified in~\cite{bombieri} and are listed in Table~\ref{classificationHyperelliptic}.

\begin{table}[H]
\begin{tabular}{@{}lllll@{}}
\toprule
 G \hspace{1,5cm}	& Lattice  \hspace{1cm}	& Generators of the action		 	\\ \midrule
 $\Z_2$			& arbitrary  		& $(z_1,z_2)\mapsto(-z_1,z_2+1/2)$ 		\\
 $\Z_2\times\Z_2$	& arbitrary 		& $(z_1,z_2)\mapsto(-z_1,z_2+1/2)$	 	\\
			&			& $(z_1,z_2)\mapsto(z_1+1/2,z_2+i\beta_2/2)$	\\
 $\Z_3$			& $\Z\oplus\omega\Z$	& $(z_1,z_2)=(\omega z_1,z_2+1/3)$		\\
 $\Z_3\times\Z_2$	& $\Z\oplus\omega\Z$	& $(z_1,z_2)=(\omega z_1,z_2+1/3)$		\\
			&			& $(z_1,z_2)=(z_1+1/2,z_2+i\beta_2/2)$		\\
 $\Z_4$			& $\Z\oplus i\Z$	& $(z_1,z_2)=(i z_1,z_2+1/4)$			\\
 $\Z_4\times\Z_2$	& $\Z\oplus i\Z$	& $(z_1,z_2)=(i z_1,z_2+1/4)$			\\
			&			& $(z_1,z_2)=(z_1+1/2,z_2+i\beta_2/2)$		\\
 $\Z_6$			& $\Z\oplus\omega\Z$	& $(z_1,z_2)=(-\omega z_1,z_2+1/6)$		\\
\bottomrule
\end{tabular}
\caption{\label{classificationHyperelliptic}List of hyperelliptic surfaces from \cite{bombieri}.}
\end{table}

We see that $G$ acts on $T^2_\Lambda$ either by rotations of multiples of the angle $\theta_k=2\pi/k$, where $k$ is the order of the first factor of $G$ (i.e. $\theta_k$ can take the values $\pi,\,3\pi/2,\,\pi/4,$ and $-\pi/3$), or by translation $x_1\mapsto x_1+1/2$.
Let $n_k=2k+2$ so that
\begin{equation}\label{EQN_CondOnTheta}
 \frac{n_k\pi}{k}=\theta_k+2\pi.
\end{equation}
We define
\begin{align*}
  X&=\del_{x_2}-\sin(n_k\pi x_2)\ \del_{x_1}+\cos(n_k\pi x_2)\ \del_{y_1}=\del_{x_2}+R(n_k\pi x_2)\del_{y_1}\\
  JX&=\del_{y_2}+\cos(n_k\pi x_2)\ \del_{x_1}+\sin(n_k\pi x_2)\ \del_{y_1}=\del_{y_2}+R(n_k\pi x_2)\del_{x_1},  
\end{align*}
where $R(\eta)$ denotes the rotation matrix of an angle $\eta$ in the plane $(x_1,y_1)$.
One can verify that this defines a \mbox{$J$-Engel} structure on $\C^2$ (see Section~\ref{SEC_Tori}). Moreover it passes to the quotient $T^2_\Lambda\times T^2_{\tilde\Lambda}$ and, being invariant with respect to the action of $G$, it defines a \JE structure on $M$. Indeed the tangent map to the first generator $g\in G$ acts on $X$ as follows
\begin{align*}
 (T_pg)(X(p))&=\del_{x_2}+R({\theta_k})R({n_k\pi x_2})\del_{y_1}=\del_{x_2}+R({\theta_k})R({n_k\pi x_2})\del_{y_1}\\
 &=\del_{x_2}+R({n_k\pi x_2+\theta_k})\del_{y_1}\overset{\eqref{EQN_CondOnTheta}}{=}\del_{x_2}+R({n_k\pi (x_2+1/k)})\del_{y_1}=X(g(p)),
\end{align*}
while $X$ is always invariant with respect to the action of the second generator of $G$.

Alternatively one can appeal to Theorem~\ref{THM_hasegawa} and construct a left-invariant \mbox{$J$-Engel} structure on the solvable group. An example is given by $\sD=\bra A ,\, JA\ket$ where $A=X_1+X_4$.
It is easy to check that
\begin{align*}
[A,\, JA]&=X_1\, \,\,\text{and} \\
[A,\, [A,\, JA]]&=-X_2\, ,
\end{align*}
so that $\sD$ is in fact a \JE structure.
Let $\{a_1,\,a_2,\,a_3,\,a_4\}$ denote the dual basis on the dual Lie algebra, then $\alpha=a_2+a_3$ and $\beta=\alpha\circ J=a_1-a_4$ are K-Engel forms and $R=X_3$ is a \JE vector field.

\subsection{Kodaira surfaces}\label{SEC_Kodaira} 
We can produce explicit \JE structures on Kodaira surfaces making use of their structure of solvmanifold.
Let us consider first primary Kodaira surfaces, these are quotients of $Nil_3\times\R$ by a cocompact lattice.
There are no geometric Engel structures on these manifolds (see Section 3 in \cite{vogelGeometricEngel}), nonetheless we can give some explicit examples of $J$-Engel structures.
Let $t$ be the coordinate on the second factor $\R$ and consider the complex plane field $
  \sD=\bra A ,\, JA\ket$ where $A=X_4+\sin t\,X_1-\cos t\,X_2$.
In this case one gets 
\begin{align*}
[A,\, JA]&=X_3-\sin t\,X_1+\cos t\,X_2\,\, \text{and} \\
[A,\, [A,\, JA]]&=-\cos t\,X_1-\sin t\,X_2\, .
\end{align*}
A straightforward computation shows that these two vectors, together with $A$ and $JA$, span the Lie algebra, proving that $\sD$ is a \JE structure.

Now for secondary Kodaira surfaces consider the Lie algebra from Theorem~\ref{THM_hasegawa} and set $\sD=\bra A ,\, JA\ket$ where $A=X_1+X_4$.
Here we see that 
\begin{align*}
[A,\, JA]&=-X_3+X_1\,\, \text{and} \\
[A,\, [A,\, JA]]&=-X_2\, ,
\end{align*}
which implies that $\sD$ is a left-invariant \JE structure.

\subsection{Complex tori}\label{SEC_Tori}

 We denote by $T_\Lambda^4$ the complex torus obtained by quotienting $\C^2$ via the lattice $\Lambda$ generated by the vectors
 \[
  e_1=\Matrix{1\\0\\0\\0},\qquad e_i=\Matrix{\alpha_i\\\beta_i\\\gamma_i\\\delta_i}
 \]
 with $\alpha_i,\,\beta_i,\,\gamma_i,\,\delta_i\in\R$ for $i=2,3,4$.
 All left-invariant distributions on $\R^4$ are integrable, so there can be no geometric Engel structure.
 
  Suppose that $\Lambda$ is such that $\alpha_i\in\Q$ for all $i$, write $\alpha_i=p_i/q_i$.
    Define the function $\theta\colon \C^2\longrightarrow \R$ via
  \[
     \theta(x_1,y_1,x_2,y_2)\colon =2\pi q_1q_2q_3x_1.
  \]
  By the hypothesis on $\Lambda$ this function passes to the quotient $T^4_\Lambda$.
  Now consider the vector fields
  \[
   A=\del_{x_1}+\sin\theta\ \del_{x_2}-\cos\theta\ \del_{y_2},\qquad JA=\del_{y_1}+\cos\theta\ \del_{x_2}+\sin\theta\ \del_{y_2}\, .
  \]
  By construction these are invariant under the action of $\Lambda$ so that they pass to $T^4_\Lambda$.
  A direct calculation shows that $\sD=\bra A,\,JA\ket$ is a \JE structure.

 Observe that the previous example gives an explicit \JE structure on a dense set of Abelian varieties $T^4_\Lambda$.
 Since the Engel condition is open under perturbations, we obtain a family of \JE structures on an open dense set of tori.

\subsection{Non-K\"ahler properly elliptic surfaces}\label{SEC_Kod1}

A properly elliptic surface is an elliptic fibration over a good orbifold (cf~\cite[p.~139]{wall}).
By Theorem~\ref{THM_geometric}  we have that all non-K\"ahler surfaces of Kodaira dimension 1 admit a geometric structure of type $G=\widetilde{SL}(2,\R)\times\R$, hence it suffices to produce a left-invariant \mbox{$J$-Engel} structure on this Lie group.

We can fix a basis $\{X_1,\,X_2,\,X_3,\,X_4\}$ for the Lie algebra $\fg$ of $G$ for which the only non-zero Lie brackets are
\[
 [X_1,X_2]=X_3,\qquad[X_2,X_3]=X_1,\qquad[X_3,X_1]=-X_2,
\]
and the complex structure $J$ is given by $JX_1=X_2$ and $JX_3=X_4$.

We can define a left-invariant complex plane field on $G$ by $\sD=\bra A,\,JA\ket$ with $A=X_1+X_2+X_3$.
By computing
\begin{align*}
[A,\, JA]&=X_2+2X_3-X_1\quad\text{and} \\
[A,\, [A,\, JA]]&=X_1+3X_2\, .
\end{align*}
one sees that $\sD$ defines a \JE structure.


\end{document}